\documentclass[12pt]{amsart}

\usepackage[margin=1.15in]{geometry}

\usepackage{amsmath,amscd,amssymb,amsfonts,latexsym,wasysym, mathrsfs, enumitem, mathtools,hhline,color}
\usepackage[all, cmtip]{xy}

\usepackage{url}

\definecolor{hot}{RGB}{65,105,225}

\usepackage[pagebackref=true,colorlinks=true, linkcolor=hot ,  citecolor=hot, urlcolor=hot]{hyperref}

\usepackage{graphicx}

\newcommand{\C}{\mathbb{C}}
\newcommand{\Q}{\mathbb{Q}}
\newcommand{\Z}{\mathbb{Z}}

\theoremstyle{plain}
\newtheorem{theorem}{Theorem}[section]
\newtheorem{prop}[theorem]{Proposition}

\newtheorem{lm}[theorem]{Lemma}

\newtheorem{cor}[theorem]{Corollary}

\newtheorem{thrm}[theorem]{Theorem}

\theoremstyle{definition}

\newtheorem{defn}[theorem]{Definition}

\newtheorem{rmk}[theorem]{Remark}

\newtheorem{ex}[theorem]{Example}
\newtheorem*{ex*}{Example}

\def\be{\begin{equation}}
\def\ee{\end{equation}}

\def\bt{\begin{thrm}}
\def\et{\end{thrm}}

\def\bc{\begin{cor}}
\def\ec{\end{cor}}

\def\br{\begin{rmk}}
\def\er{\end{rmk}}

\def\bp{\begin{prop}}
\def\ep{\end{prop}}

\def\bl{\begin{lm}}
\def\el{\end{lm}}

\def\bex{\begin{ex}}
\def\eex{\end{ex}}

\def\bd{\begin{defn}}
\def\ed{\end{defn}}

\newcommand\sK{{\mathcal K}}
\newcommand\sH{{\mathcal H}}
\newcommand\sC{{\mathscr C}}

\newcommand\sF{{\mathcal F}}

\newcommand\sM{{\mathcal M}}

\newcommand\sV{{\mathcal V}}

\newcommand\sR{\mathcal{R}}

\newcommand\sS{\mathcal{S}}
\def\sQ{\mathcal{Q}}

\newcommand{\mb}{\mathcal{M}_{\textrm{B}}}

\DeclareMathOperator{\Exp}{Exp}
\DeclareMathOperator{\homo}{Hom}

\DeclareMathOperator{\spec}{Spec}

\DeclareMathOperator{\Perv}{Perv}
\DeclareMathOperator{\Mod}{Mod}

\def\ra{\rightarrow}

\def\bC{\mathbb{C}}
\def\cM{\mathcal{M}}

\def\al{\alpha}

\def\cH{\mathcal{H}}

\def\cO{\mathcal{O}}
\def\lra{\longrightarrow}
\def\bQ{\mathbb{Q}}
\def\ol{\overline}
\def\cL{\mathcal{L}}

\def\bZ{\mathbb{Z}}

\def\ul{\underline}
\def\lam{\lambda}

\def\sD{\mathscr{D}}
\def\xra{\xrightarrow}
\def\bDbc{{\mathbf{D}}^b_c}
\def\bDrh{{\mathbf{D}}^b_{rh}}
\def\cA{\mathcal{A}}
\def\cF{\mathcal{F}}
\def\cS{\mathcal{S}}
\def\rb{\mathcal{R}_B}

\title{On the length of perverse sheaves and $\sD$-modules}

\author{Nero Budur}
\address{KU Leuven, 
Celestijnenlaan 200B, B-3001 Leuven, Belgium}
\email{nero.budur@kuleuven.be}

\author{Pietro Gatti}
\address{KU Leuven, 
Celestijnenlaan 200B, B-3001 Leuven, Belgium, and University of Padua, Torre Archimede, Via Trieste 63, 35121 Padova, Italy}
\email{pietro.gatti@kuleuven.be} 

\author{Yongqiang Liu}
\address{KU Leuven, 
Celestijnenlaan 200B, B-3001 Leuven, Belgium. Currently: BCAM, Mazarredo 14, 48009 Bilbao Basque Country, Spain}
\email{yliu@bcamath.org}

\author{Botong Wang}
\address{University of Wisconsin, Van Vleck Hall, 480 Lincoln Drive, Madison, WI, USA} 
\email{bwang274@wisc.edu}

\begin{document}

\maketitle

\begin{abstract} We prove that the length function for perverse sheaves and algebraic regular holonomic $\sD$-modules on a smooth complex algebraic variety $Y$ is an absolute $\bQ$-constructible function. One consequence is: for ``any" fixed natural (derived) functor $F$ between constructible complexes or perverse sheaves on two smooth varieties $X$ and $Y$,  the loci of rank one local systems $L$ on $X$ whose image $F(L)$ has prescribed length are Zariski constructible subsets defined over $\bQ$, obtained from finitely many torsion-translated complex affine algebraic subtori of the moduli of rank one local systems via a finite sequence of taking union, intersection, and complement.
 \end{abstract}

\section{Introduction}

Let $X$ be a complex algebraic variety.  It is well known that the category of perverse sheaves $\Perv(X,\C)$ on $X$ is artinian and noetherian. In other words, every perverse sheaf $P$ has a finite length filtration, called a composition series, $$ 0=P_{0}\hookrightarrow P_{1} \hookrightarrow \cdots \hookrightarrow P_{\ell}=P $$
such that the quotients $P_{i}/P_{i-1}$  are non-zero simple perverse sheaves for $i=1,\ldots, \ell$. Moreover, every simple perverse sheaf on $X$ is an intersection complex $IC_{\overline{S}}(L)$ of some   irreducible local system $L$ on a smooth irreducible subvariety $S$ of $X$, where $\overline{S}$ denotes the closure of $S$. The length of $P$ is denoted by $\ell(P)$ and is by definition  the number $\ell$ of simple constituents in a composition series. 

If $X$ is in addition smooth, let $\Mod_{rh}(\sD_X)$ denote the category of algebraic regular holonomic  $\sD_X$-modules on $X$. This category is also artinian and noetherian, and hence it admits a length function. By the Riemann-Hilbert correspondence, there is an equivalence of categories 
$$
DR_X:\Mod_{rh}(\sD_X) \xra{\sim} \Perv(X,\C),
$$
hence the length of a regular holonomic $\sD$-module agrees with that of its perverse sheaf counterpart.  

For an introduction to perverse sheaves, $\sD$-modules, and the Riemann-Hilbert correspondence, one has now  as classical references the books \cite{D}, \cite{HTT}; see also \cite{DM}.

Despite the two different points of view, computing lengths of perverse sheaves or $\sD$-modules remains a difficult task. There is currently no available algorithm.  A basic open question is for example to compute the length of $\sD_X$-modules generated by complex powers $f^\alpha$ for some polynomial $f:X=\bC^n\ra\bC$ and $\al\in\bC$, see U. Walther \cite{Wa}. Answers are now available for $\sD_Xf^{-1}=\cO_X[1/f]$, this case being also related to local cohomology, for $f$ a hyperplane arrangement by Abebaw-B\o gvad \cite{AB},  \`Alvarez Montaner - Garc\'ia L\'opez - Zarzuela Armengou \cite{AGZ}, B\o gvad-Gon\c{c}alves \cite{BG}, Budur-Saito \cite{BS}, T. Oaku \cite{O}, U. Walther \cite{W2}; and also for any $\al$ for $\sD_Xf^\al$ in the case $f$ is a quasi-homogeneous polynomial with an isolated singularity by Bitoun-Schedler \cite{BiS}.
In terms of perverse sheaves, these $\sD$-modules correspond to derived push-forwards of rank one local systems from the complement of the zero locus of $f$, see equation (\ref{eqER}) below. We will adopt the topological language for the rest of the introduction.

The purpose of this article is to show that  collections of objects with prescribed length are very special. The most concrete versions of this statement start with the simplest kind of perverse sheaves, the $\bC$-local systems of rank one (or, equivalently for flat line bundles if $X$ is smooth), since the moduli space of rank one local systems 
$$
\mb(X,1)=\homo(H_1(X,\bZ),\bC^*)
$$
is an easy space to grasp. It is a complex affine torus, that is, an algebraic group isomorphic to the product of some multiplicative group $(\bC^*)^{b}$ with a finite abelian group. 

\begin{defn}\label{defAr1} An {\it absolute $\bQ$-constructible subset} of $\mb(X,1)$ is a Zariski constructible subset defined over $\bQ$ which is obtained from finitely many torsion-translated complex affine algebraic subtori  via a sequence of taking union, intersection, and complement.
\end{defn}

Here is one sample from the results in this article:

\bt\label{thrmMain} Let $f: X \to Y$ be an algebraic map between two complex algebraic varieties. Let $S$ be an irreducible smooth subvariety of $X$. Then for any $i,k\in\bZ$, the sets
$$\{ L \in \sM_B(S, 1)  \mid  \ell( {}^pR^if_*(IC_{\ol{S}}(L)) ) = k \},$$
$$\{ L \in \sM_B(S, 1)  \mid  \ell(  {}^pR^if_!(IC_{\ol{S}}(L))  ) = k \},$$
are absolute $\bQ$-constructible,  where ${}^pR^if_*={}^p\cH^i\circ Rf_*$ and $ ^pR^if_!={}^p\cH^i\circ Rf_!$.
\et

Here ${}^pR^if_*$ and ${}^pR^if_!$ are the functors called perverse $i$-th derived direct image and, respectively, direct image with compact supports. They are defined on the  bounded derived category of constructible sheaves of $\bC$-vector spaces on $X$, denoted
$$
\bDbc (X,\bC).
$$

With same assumptions as Theorem \ref{thrmMain}, if $f$ is in addition proper then  the  Decomposition Theorem  says that ${}^pR^if_*(IC_{\ol{S}}(L))={}^pR^if_!(IC_{\ol{S}}(L)) $ is a semi-simple perverse sheaf. Hence in this case 
the length counts exactly the number of direct summand factors.  
  
The definition of length can be extended to the whole derived category using the  perverse cohomology sheaf functors $${}^p\cH^i:\bDbc(X,\bQ)\ra \Perv(X,\bC).$$ For $\sF$ in  $\bDbc(X,\C)$ one can define
$$ \ell(\sF)= \sum_i (-1)^i \ell (^p \sH^i (\sF)),$$
since  ${}^p\cH^i(\sF)$ is non-zero for only finitely many $i$.

The length function is additive for distinguished triangles in $\bDbc(X,\bC)$. If $X$ is a point, then $\ell(\sF)$ is just the classical Euler characteristic for a bounded complex of  $\C$-vector spaces with finite-dimensional cohomologies.
Since $${}^pR^if_*={}^p\cH^i\circ Rf_*\quad\text{ and }\quad ^pR^if_!={}^p\cH^i\circ Rf_!$$ by definition, Theorem \ref{thrmMain} also holds true with ${}^pR^if_*$ and ${}^pR^if_!$ replaced by $Rf_*$ and $Rf_!$, respectively.

Under an affine open embedding, the derived direct image (with compact supports) functor is perverse $t$-exact. A direct implication of the above theorem is the following  basic yet interesting case, obtained by embedding a possibly singular variety $X$ into a larger smooth variety:

\begin{cor}\label{corJ}  Let $j : U \to X$ be an affine open embedding (e.g. the complement of a hypersurface) of complex algebraic
varieties with $U$ smooth of dimension $n$. Then for any $i,k\in\bZ$, the set
$$\ell_k(U,X):=\{ L \in \sM_B(U, 1)  \mid  \ell(Rj_*(L[n])) = k \}$$
is absolute $\Q$-constructible.
\end{cor}

When $U=X\setminus\{f=0\}$ for a regular function $f:X\ra\bC$, and $\al\in\bC$, then the complex-powers $\sD$-module $\sD_X f^\al$ corresponds under the Riemann-Hilbert equivalence to $Rj_*(L[n])$ with $L$ the kernel of the flat connection on the $\cO_U$-invertible sheaf $\sD_uf^\al$. Hence in this case,
\begin{equation}\label{eqER}
\ell(Rj_*(L[n]))=\ell (\sD_X f^\al).
\end{equation}

\begin{rmk}
The locus $\ell_1(U,X)$ of rank one local systems on $U$ with simple maximal perverse extension to $X$ has been addressed in \cite{BLSW} and some deeper properties have been found. Below we assume that $U$ is the complement of a hypersurface in a complex manifold $X$. Then:
\begin{itemize}
\item $\ell_1(U,X)$ is {Zariski} open.
\item Denote the complement of $\ell_1(U,X)$ in $\mb(U,1)$ by $\ell_{\ge 2}(U,X)$. (This was denoted $\sV^{ns}(U,X)$ in \cite{BLSW}.) Then $\ell_{\ge 2}(U,X)$ can be expressed in terms of cohomology support loci of rank one local systems on small ball complements along the hypersurface $X\setminus U$, \cite[Theorem 1.5]{BLSW}. These are known to be finite unions of torsion-translated complex affine algebraic subtori even in this analytic setup by \cite{BW16}.
\item If $X\setminus U=\cup_jf_j^{-1}(0)$ for some non-invertible analytic functions $f_j$ on $X$ with $j=1,\ldots, r$, denote by $\sS(F)$ the inverse image  of $\ell_{\ge 2}(U,X)$ in $\mb((\bC^*)^r)\simeq (\bC^*)^r$, under the map pulling back local systems from  $(\bC^*)^r$ to $U$  via 
$$
F=(f_1,\ldots,f_r):U\ra (\bC^*)^r.
$$
Then $\sS(F)$ is pure of codimension one in $(\bC^*)^r$, and each irreducible component of $\sS(F)$ is a torsion-translated complex affine algebraic subtorus, \cite[Theorem 1.3 and Corollary 1.6]{BLSW}.
\item It is not true that $\ell_{\ge 2}(U,X)$ is always of pure codimension one, \cite[Example 3.2]{BLSW}.
\item $\sS(F)$ can be computed from any log resolution of $\prod_jf_j$ \cite[Theorem 1.3]{BLSW}.
\item $\ell_{\ge 2}(U,X)=\sS(F)$ in the case that each $f_j$ is linear, that is $U$ is the complement of a hyperplane arrangement, and there is a simple combinatorial formula for $\ell_{\ge 2}(U,X)$, \cite[Theorem 1.2]{BLSW}. 
\item Conjecturally, $$\sS(F)=\Exp (Z(B_F)),$$ where $B_F\subset \bC[s_1,\ldots,s_r]$ is the Bernstein-Sato ideal of $F$, $Z(B_F)$ denotes its zero locus in $\bC^r$, and $\Exp$ is the map $\exp (2\pi i (\_))$ coordinate-wise. There are algorithms already implemented for computing Bernstein-Sato ideals. Only the inclusion $\sS(F)\subset \Exp(Z(B_F))$ is currently known, \cite{B-ls}. 
\item The conjectured equality $\sS(F)=\Exp (Z(B_F))$ implies a transfer of the above properties of $\sS(F)$ onto the Bernstein-Sato ideal: the codimension-one irreducible components of the zero locus of $B_F$ must be linear, their union must be defined over $\bQ$, and each component of codimension $\ge 2$  must be up to translation by a vector in $\bZ^r$ included in the codimension-one locus. These conjectural properties of $B_F$ have recently been confirmed by Maisonobe \cite{Ma}.
\end{itemize}
\end{rmk}

It remains open whether some of these properties can be extended to the higher-length loci $\ell_k(U,X)$. Currently we cannot prove or disprove that $\ell_{\ge k}(U,X)$ is closed for $k\ge 3$, for example. For an explicit $\sD$-module theoretic interpretation of the sets $\ell_k(U,X)$, see \cite{BLSW}.

Theorem \ref{thrmMain} is only a sample illustration of a very general result. This result relies on the theory of absolute constructible sets of \cite{BW17}. Firstly, in \cite{BW17} it was defined more generally what a $\bQ$-constructible and an absolute $\bQ$-constructible set of objects up to isomorphisms in $\bDbc(X,\bC)$ and $\Perv(X,\bC)$ is for a complex algebraic variety $X$. This general notion of constructibility frees the classical notion of constructibility from needing to live on a moduli space parametrizing a portion of the category, but agrees with the classical notion when a moduli space exists. The notion of absoluteness involves the full Riemann-Hilbert correspondence with the derived category of algebraic regular holonomic $\sD$-modules, and therefore needs the extra assumption that $X$ is smooth.  A difficult result was to characterize all the absolute $\Q$-constructible sets of $\mb(X,1)$ as in Definition \ref{defAr1}. Work in progress is dedicated now to the characterization of such sets for  moduli of other small-rank local systems. 

In \cite{BW17} it was shown that ``any natural" functor between categories such as $\bDbc(X,\bC)$ and $\Perv(X,\bC)$ for smooth algebraic varieties is absolute $\bQ$-constructible. That is, it has the property that it pulls back an absolute $\bQ$-constructible set to another one. 
Smoothness may be bypassed occasionally at intermediate steps by embedding a variety $X$ into a bigger smooth variety $Y$ and considering constructible complexes and perverse sheaves on $Y$ with support on $X$. In practice, one could not define what ``any natural" functor means, and \cite{BW17} resorted to providing a huge list of functors in the hope that via usual constructions they might generate any other desired natural functor. One of the functors that was missing from that list is the length function. In this article we add the length function to the list of absolute  $\bQ$-constructible functions for smooth complex algebraic varieties:

\bt \label{mainAbs}  Let $X$ be a complex algebraic variety. Then:
\begin{itemize}
\item[(a)]
the function $\ell:Perv(X, \C) \ra \bZ$ is  $\Q$-constructible.
\item[(b)]
the function $\ell:\bDbc(X, \C)  \ra \Z $ is  $\Q$-constructible. 
\end{itemize}
If moreover $X$ is smooth, then both the functions  are absolute  $\Q$-constructible.
\et

The proper definitions will be recalled in the next section. 

This theorem implies that any mixture of functors as in \cite[Theorem 5.14.1]{BW17} on complex algebraic varieties, now including also the length function,  is a $\bQ$-constructible functor, or a $\bQ$-constructible function if it ends in $\bZ$. Further, by \cite[Theorem 6.4.3]{BW17}, if all the functors in the mixture come from smooth varieties, then the resulting mixture functor or function is also absolute. By mixture, we mean usual operations one can do with functors, such as composition and fiber products. In particular:

\begin{thrm}\label{thrmMainG} Let $X$ be a smooth complex algebraic variety. Let $F$ be a map from $\mb(X,1)$ to  constructible complexes or perverse sheaves up to quasi-isomorphisms on another smooth variety $Y$.  Assume that the map $F$ is as in \cite[Theorem 6.4.3]{BW17}. Then the subset of $\mb(X,1)$ given by
$$
\{L\in\mb(X,1)\mid \ell(F(L))=k\}
$$
is absolute $\bQ$-constructible.
 \end{thrm} 

Theorem \ref{thrmMain} follows by setting $F={}^p\cH^i\circ Rf_*$ and $F={}^p\cH^i\circ Rf_!$.  

The proof of Theorem \ref{mainAbs} boils down to \cite[Theorem 6.4.3]{BW17} by using induction on the dimensions of the strata pertaining to a constructible complex. This is done in Section 3. In Section 2 we recall the basic material on absolute constructibility from \cite{BW17}. In the last section we provide some examples of length jump loci of rank 2 local systems. 

\smallskip
 {\it Acknowledgement.} We thank L. Saumell and J. Sch\"urmann for many useful discussions. The first three authors were partly sponsored by the research grant STRT/1 3/005 and a Methusalem grant METH/15/026 from KU Leuven,  and by the research projects G0B2115N, G0F4216N, G097819N from the Research Foundation of Flanders.  The third author is partly supported by the ERCEA 615655 NMST Consolidator Grant, the Basque Government through the BERC 2018-2021 program, and by the Spanish Ministry of Science, Innovation and Universities: BCAM Severo Ochoa accreditation SEV-2017-0718. The first author thanks J. Fernandez de Bobadilla and BCAM Bilbao for the hospitality during writing part of this article via the fund ERC-615655-NMST.

\section{Absolute $\bQ$-constructibility}

In this section we summarize some of the theory of absolute constructibility from \cite{BW17} that we  use: 

\subsection{A general notion of constructibility: unispaces.} We recall first the general notion of constructibility defined in \cite{BW17}. The idea is
that in order to define constructibility one can use a much more naive notion of a family of objects than the usual flatness condition. In fact only a generic base change property should play the essential role, instead of a full base change.

Let $\cA lg_{ft,reg}(\bC)$ be the category of finite type regular $\bC$-algebras. Let $R\in\cA lg_{ft,reg}(\bC)$.  Let $X$ be a complex algebraic variety and $Y$ a $\bC$-scheme of finite type.
We will be mainly concerned with:
\be\label{eqSp}
Y(R), \bDbc(X,R),\text{ and } \Perv(X,R).
\ee 
Here, $Y(R)$ denotes the set of $R$-points of $Y$, that is, morphisms $\spec(R)\ra Y$ over $\bC$. Also, $\bDbc(X,R)$ denotes the triangulated full subcategory of the derived category of sheaves of $R$-modules on the analytic variety underlying  $X$ consisting of complexes with finitely many non-zero cohomology sheaves and such that the stalks of the cohomology sheaves are finitely generated $R$-modules. Further, using a non-standard definition, we let
\begin{equation}\label{eqP}
\Perv(X,R)=\{ \cF_R\in\bDbc(X,R)\mid \cF_R\otimes^L_RR/m\in\Perv(X,\bC)\text{ for all }m\in\spec(R)(\bC)\},
\end{equation}
where a $\bC$-point $m$ of $\spec(R)$ is identified with the associated maximal ideal of $R$.
The (set of isomorphism classes of objects in) the three types of ``spaces" in (\ref{eqSp}) are endowed with a natural base change maps under a morphism $R\ra R'$ in $\cA lg_{ft, reg}(\bC)$:
$$
Y(R)\ra Y(R'),$$ 
$$\bDbc(X,R)_{/\simeq}\xra{\otimes^L_RR'}\bDbc(X,R')_{/\simeq},$$
$$ \Perv(X,R)_{/\simeq}\xra{\otimes^L_RR'}\Perv(X,R')_{/\simeq}.
$$
These are typical examples {\it unispaces} in the terminology of \cite{BW17}, where the {\it category of unispaces} was defined. A unispace $(\sC(\_),\star)$ is just a functor
$$
\sC:\cA lg_{ft,reg}(\bC)\lra \cS et
$$
$$
R\mapsto \sC(R)
$$
$$
(R\ra R')\mapsto (\sC(R)\xra{\star_RR'}\sC(R')).
$$
A {\it morphism of unispaces} $F:\sC\ra\sC'$ is a map $F(\bC):\sC(\bC)\ra\sC'(\bC)$ with certain compatibility properties, \cite[2.4.5]{BW17}. We skip the precise definition but we recall the most important properties.

\begin{defn} Let $(\sC(\_),\star)$ be a unispace. A function
$$
\phi:\sC(\bC)\lra \bZ
$$
is a {\it constructible function with respect to $\sC$} if 
$$
\forall\; R\in \cA lg_{ft,reg}(\bC)\text{ and }\forall\; \cF_R\in \sC(R)$$
the induced function on the $\bC$-points
$$
\spec(R)(\bC)\ra\bZ, \quad m\mapsto\phi(\cF_R\star_RR/m)
$$
is constructible in the classical sense.   
\end{defn}

The category of unispaces is tailor-made to form the most natural environment for constructible functions:

\begin{thrm}{\rm{(}\cite{BW17}\rm{)}}
\begin{enumerate}
\item If $Y$ is a $\bC$-scheme of finite type (or an algebraic stack), then a function
$
Y(\bC)\ra \bZ
$ is constructible in the classical sense if and only if it is constructible with respect to the unispace $Y(\_)$.
\item A function $\sC(\bC)\ra \bZ$  constructible with respect to a unispace $\sC$ is the same as a morphism of unispaces $\sC\ra\ul{\bZ}^{cstr}$, where the unispace $\ul{\bZ}^{cstr}$ is defined by setting
$
\ul{\bZ}^{cstr}(R)
$ to be the set of all constructible functions $\spec(R)(\bC)\ra\bZ$, with the natural base change maps $\ul{\bZ}^{cstr}(R)\ra\ul{\bZ}^{cstr}(R')$ for  morphisms $R\ra R'$ in $\cA lg_{ft,reg}(\bC)$.
\item Let $\sC$ and $\sC'$ be two unispaces. Let $F(R):\sC(R)\ra\sC'(R)$ be a collection of maps, one for every $R\in\cA lg_{ft,reg}(\bC)$. Assume that for all integral domains $R\in\cA lg_{ft,reg}(\bC)$ and every $\cF_R\in\sC(R)$, the base change formula
$$
(F(R)(\cF_R))\star_RR/m=F(R/m)(\cF_R\star_RR/m)
$$
of elements in $\sC'(R/m)=\sC'(\bC)$ holds for all $m$ in a non-empty open dense subset of $\spec(R)(\bC)$. Then $F(\bC)$ lifts uniquely to a morphism of unispaces $F:\sC\ra\sC'$. 
\end{enumerate}
\end{thrm}

Part (3) above is how one constructs morphisms of unispaces in practice. For example, given a functor $\bDbc(X,\bC)\ra\bDbc(Y,\bC)$, it usually admits a version with $R$-coefficients, and if it satisfies the generic base change of (3), one gets a morphism of unispaces.
Since the unispaces form a category, composition of two morphisms is a morphism, and so morphisms of unispaces preserve constructibility of functions by part (2) above. Then one produces classical constructible functions by applying part (1). For example, to produce classical constructible sets on moduli spaces $\mb(X,r)$ of local systems  of rank $r$, one composes the associated morphism of unispaces  $\rb(X,r)(\_)\ra\bDbc(X,\_)$  with a constructible function on the latter unispace, where $\rb(X,r)$ is the fine moduli space of $r$-dimensional representations of the fundamental group of $X$. Then one uses the morphism of $\bC$-schemes of finite type $\rb(X,r)\ra\mb(X,r)$ to push down to $\mb(X,r)$ such classical constructible sets constructed by (derived) functors.

In defining the new notion of constructibility in all the above we worked with finite type regular $\bC$-algebras. By working with finite type regular $\bQ$-algebras, one can define the category of unispaces defined over $\bQ$, and use it to produce in a non-classical way classical $\bQ$-constructible subsets of varieties defined over $\bQ$, \cite[Section 3]{BW17}. For example $Rf_*$ gives a morphism defined over $\bQ$ of unispaces defined over $\bQ$ since $Rf_*$ exists already on $\bDbc(X,\bQ)$. Similarly the functor $(\_)\otimes^L_\bC \cF$  gives a morphism defined over $\bQ$ of unispaces defined over $\bQ$ if $\cF$ is  in $\bDbc(X,\bQ)$, as opposed to $\cF$ being merely an object in $\bDbc(X,\bC)$.

Theorem 5.14.1 of \cite{BW17} gives a  list of functors, including all the basic ones which first come to one's mind, which give morphisms of unispaces and shows that most of them are defined over $\bQ$. We shall refer to that list in the course of the proof of the main result.

\subsection{Absoluteness.} Let $X$ be now a smooth complex algebraic variety. The full Riemann-Hilbert correspondence states that there are equivalences of categories
$$
\xymatrix{
\cM od_{rh}(\sD_X) \ar[d] \ar[r]^\sim_{RH}  & \Perv(X,\bC) \ar[d]\\
\bDrh(\sD_X)  \ar[r]^{\sim}_{RH}  & \bDbc(X,\bC)
}
$$
compatible with the usual functors, where the left-hand side is the (bounded derived) category of algebraic regular holonomic $\sD_X$-modules, see \cite{HTT}, and the vertical arrows are embeddings of full subcategories. We will use algebraicity of the $\sD$-modules to conjugate them by any element $\sigma$ of $Gal(\bC/\bQ)$. The conjugate of a $\sD_X$ module is a $\sD_{X^\sigma}$-module, where $X^\sigma$ is the conjugate of $X$.

\begin{defn}
A function $\phi:\bDbc(X,\bC)_{/\simeq}\ra\bZ$ that is $\bQ$-constructible with respect to the unispace $\bDbc(X,\_)$ is called {\it absolute} if for all $\sigma\in Gal(\bC/\bQ)$, the function $\phi^\sigma:\bDbc(X^\sigma,\bC)_{/\simeq}\ra\bZ$, associated to $\phi$ via  conjugation by $\sigma$ of the algebraic $\sD$-modules, is also $\bQ$-constructible with respect to the unispace $\bDbc(X^\sigma,\_)$.
\end{defn}

Similarly, one can replace in this definition $\bDbc(X,\bC)$ with $\Perv(X,\bC)$, or with the category of $\bC$-local systems on $X$, etc., to get a notion of absoluteness in those contexts. 

Theorem 6.4.3 of \cite{BW17} gives a  list of functors, including all the basic ones, which give morphisms of unispaces defined over $\bQ$ that preserve  absolute $\bQ$-constructible functions. Such a functor is called an {\it absolute $\bQ$-constructible functor}. Hence absolute $\bQ$-constructible functors can be used to produce absolute $\bQ$-constructible subsets of moduli spaces parametrizing portions of $\bDbc(X,\bC)$. In the case of the moduli $\mb(X,1)$ of rank one local systems, this definition of absolute $\bQ$-constructible subsets is equivalent with that of Definition \ref{defAr1}, by \cite[Theorem 8.1.2]{BW17}.

\section{Proof of Theorem \ref{mainAbs}}

Note that the complete statement of this theorem should include a reference to the unispace structures with respect to which these functions should be $\bQ$-constructible. We take these unispaces to be $\Perv(X,\_)$ and $\bDbc(X,\_)$, respectively, as defined in the previous section.

First, we claim that the parts (a) and (b) are equivalent. Suppose we know that the length function on $\bDbc(X,\bC)$ is (absolute) $\bQ$-constructible with respect to the unispace $\bDbc(X,\_)$. Consider the functor which considers perverse sheaves as objects of the derived category:
$$
\iota:\Perv(X,\bC)\lra \bDbc(X,\bC).
$$
Then the length function on $\Perv(X,\bC)$ is the composition of $\iota$ with  $\ell:\bDbc(X,\bC)\ra \bZ$. By \cite[Theorem 5.14.1]{BW17}, $\iota$ lifts to a morphism defined over $\bQ$ of unispaces defined over $\bQ$
$$\iota(\_):\Perv(X,\_)\ra\bDbc(X,\_).$$
Hence $\iota$ is a $\bQ$-constructible functor with respect to the two unispaces, and thus $\ell\circ\iota$ is a $\bQ$-constructible function on $\Perv(X,\bC)$ with respect to the unispace $\Perv(X,\_)$. By Theorem \cite[6.4.3]{BW17}, $\iota$ is also an absolute $\bQ$-constructible functor if $X$ is addition smooth, and this implies that $\ell\circ\iota$ is an absolute $\bQ$-constructible function in this case. 

Conversely, suppose  we know that $\ell:\Perv(X,\bC)\ra\bZ$ is  $\bQ$-constructible (respectivey, absolute $\bQ$-constructible if $X$ is smooth) with respect to the unispace $\Perv(X,\_)$. For a finite set $I\subset\bZ$, define 
$$
\ell_I:\bDbc(X,\bC)\ra \bZ,\quad \ell _I(\cF)=\sum_{i\in I}(-1)^i\ell({}^p\cH^i(\cF)).
$$
Then $\ell_I$ is the composition of 
$$
\bDbc(X,\bC) \xra{({}^p\cH^i)_{i\in I}} \Perv(X,\bC)^{\times I} \xra{\ell^{\times I}}\bZ^{\times I}\xra{\chi}\bZ,
$$
where $$\chi((n_i)_{i\in I})=\sum_{i\in I}(-1)^in_i.$$
Each of the functors ${}^p\cH^i$, and hence also their  product $({}^p\cH^i)_{i\in I}$,
 is an (absolute) $\bQ$-constructible functor, by \cite[Theorem 6.4.3]{BW17}. Similarly, $\ell^{\times I}$ is the product of (absolute) $\bQ$-constructible function, hence it is itself an (absolute) $\bQ$-function. It follows that $\ell_I$ is an (absolute) $\bQ$-constructible function. 
 
 Now, to show that $\ell:\bDbc(X,\bC)\ra\bZ$ is an (absolute) $\bQ$-constructible function, take $R$ a finite type regular $\bC$-algebra and $\cF_R\in\bDbc(X,R)$. We have to show that for every $k$ the set of maximal ideals
 \be\label{eqL}
\ell(\cF_R,k)= \{m\in\spec(R)(\bC)\mid \ell(\cF_R\otimes^L_RR/m)=k\}
 \ee
is constructible in the classical sense. For absoluteness, we need to show the same property similarly for $R$ replaced by a finite type regular $\bQ$-algebra, and that $\ell$ satisfies the Galois property from the definition of absoluteness, if in addition $X$ is smooth. It is clearly enough to assume that $R$ is an integral domain.

Firstly, there are only finitely many $i$ such that ${}^p\cH^i(\cF_R)$ is non-zero. Here the perverse cohomology sheaf is defined as usual, by defining first inductively the perverse cohomology truncation functors. Since we work on a complex algebraic variety, the inductive process terminates. Secondly, we know from \cite[Proposition 5.9.1]{BW17} that 
$$
{}^p\cH^i(\cF_R\otimes^L_RR/m)={}^p\cH^i(\cF_R)\otimes^L_RR/m,
$$
for $m$ outside a closed subscheme $\spec(R_1)(\bC)\subsetneq\spec(R)(\bC)$ (defined over $\bQ$) for a surjective map of $\bC$-algebras $R\ra R_1$ (defined over $\bQ$), this being the property responsible for the functor ${}^p\cH^i$ to extend to a morphisms of unispaces (defined over $\bQ$).
Hence for $m\not\in\spec(R_1)(\bC)$, there are only finitely many $i$ with ${}^p\cH^i(\cF_R\otimes^L_RR/m)\ne 0$.

By noetherian induction, there exists a disjoint  union decomposition into finitely many smooth affine locally closed subschemes $\spec(R)(\bC)=\cup_j\spec(R_i)(\bC)$ such that
$$
{}^p\cH^i(\cF_R\otimes^L_RR/m) = {}^p\cH^i(\cF_{R_j})\otimes^L_{R_j}R_j/m
$$
for $m\in\spec(R_j)(\bC)$, where $\cF_{R_j}=\cF_R\otimes^L_RR_j$. We conclude that there exists a finite set $I\subset \bZ$ such that ${}^p\cH^i(\cF_R\otimes^L_RR/m)=0$ for all $i\not\in I$ and for all $m\in\spec(R)(\bC)$.   Hence
$$
\ell(\cF_R,k)= \{m\in\spec(R)(\bC)\mid \ell_I(\cF_R\otimes^L_RR/m)=k\}=\ell_I(\cF_R,k).
$$
Since $\ell_I$ is a $\bQ$-constructible function, it follows that $\ell(\cF_R,k)$ is a $\bQ$-constructible set, and hence $\ell:\bDbc(X,\bC)\ra\bZ$ is a $\bQ$-constructible function. Clearly $\ell$ will be absolute if $X$ is smooth, since $\ell$ is compatible with the length function on $\bDrh(\sD_X)$.

This finishes the proof of the claim that parts (a) and (b) are equivalent.

Next, suppose that we have proved that the length function is $\bQ$-constructible. If in addition $X$ is smooth then the length function is automatically absolute too since the length function is compatible with the Riemann-Hilbert correspondence. Thus from now on we focus on $\bQ$-constructibility. We give the proof for constructibility over $\bC$. Since the proof for $\bQ$-constructibility is similar, we only point out a specific places in the proof where additional arguments are needed for $\bQ$-constructibility.
 
We will  prove the theorem now  by induction on the dimension of $X$. 

Let us remark first that a similar argument to the one above, based on the generic base change for the usual cohomology sheaves $\cH^i(\cF_R\otimes^L_RR/m)=\cH^i(\cF_R)\otimes^L_RR/m$, gives that for a a fixed complex $\cF_R$, there is a finite Whitney stratification of $X$ with respect to which the complexes $\cF_R$ and $\cF_R\otimes^L_RR/m$ for all $m\in\spec(R)(\bC)$ have simultaneously constructible cohomology sheaves.

If $\dim X=0$, then $X$ is just a point. In this case, $\ell:\bDbc(pt, \bC)\ra\bZ$ is the Euler characteristic of a complex of finite-dimensional $\bC$-vector spaces. This function is (absolute) $\bQ$-constructible with respect to the unispace $\bDbc(pt,\_)$ by \cite[Theorem 6.4.3]{BW17}.

Assume now that the theorem holds when $\dim X<n$. We show that it holds for $\dim X=n$.  

Take $\cF_R\in \Perv(X,R)$ for some finite type regular $\bC$-algebra  $R$, where $\Perv(X,R)$ is defined as in (\ref{eqP}). That is, $\cF_R\otimes^L_RR/m$ is in $\Perv(X,\bC)$ for all $m$. Take a  Whitney stratification of $X$ such that the cohomology sheaves of $\sF_R$ and of all $\cF_R\otimes^L_RR/m$ are locally constant on the strata.  Let $U\subset X$ be the union of the top dimensional strata. Then $U$ is open and smooth. Set $Z=X\setminus U$.

For simplicity we use the notation $\sF_m:=\sF_R\otimes^L_R R/m$. By assumption, $\cF_m$ is in $\Perv(X,\bC)$.

If $Z=\emptyset$, then $U=X$ is smooth and $\cF_R$ and $\sF_m$ are local systems, up to a shift by $n$, with finite rank for any maximal ideal $m\subset R$. Generic base change for cohomology sheaves gives again that there exists a finite disjoint union  decomposition into smooth affine locally closed subvarieties
\begin{center}
$\spec(R)(\bC)= \sqcup_j \spec(R_j)(\bC)$, 
\end{center}  such that $\sF_m$ has constant rank along $m \in \spec(R_j)(\bC)$ for fixed $j$. The rank function is indeed an (absolute) $\bQ$-constructible function, \cite{BW17}. Then we only need to show that the function $m\mapsto \ell(\cF_m)$ is constructible in the classical sense on each $\spec(R_j)(\bC)$ individually. Without loss of generality, we assume that $\sF_m$  has constant rank $r$ for any $m\in \spec(R)(\bC)$. And by generic base change for cohomology sheaves, we can also assume that $\cH^{-n}(\cF_R)$ is a local system of free finitely generated $R$-modules of rank $r$ such that $\cH^{-n}(\cF_R)\otimes^L_RR/m=\cF_m$ for all $m$. In other words, we can assume that $\cF_R$ is itself a shifted local system of free finitely generated $R$-modules. 

Since the length of a local system as a (shifted) perverse sheaf is the same as that as a local system, the problem in this case is reduced to showing that the length function on local systems is constructible with respect to the unispace
$\cL oc\cS ys_{free}(X,\_)$. This unispace is defined by setting $\cL oc\cS ys_{free}(X,R)$ to be the isomorphism classes of local systems of free finitely generated $R$-modules on $X$. Moreover, there is a  morphism of unispaces
$$
\sR(X,r)(\_)\ra \cL oc\cS ys_{free}(X,\_)
$$
where $$\sR(X,r)=Hom (\pi_1(X), GL_{r})$$ is the rank $r$ representation variety of the fundamental group of $U$, which is an affine scheme defined over $\bQ$. Since $\sR(X,r)$ is a fine moduli space and the length function on local systems is the same as the length of a representation, it is enough to show that the length function on the $\bC$-points of  $\sR(X,r)$ is $\bQ$-constructible in the classical sense. This is of course well-known already.

When $Z\neq \emptyset$, let $j$ and $i$ denote the inclusions from $U$ and $Z$ to $X$, respectively. By replacing $U$ with a smaller affine open subset, we can assume that the inclusion $j$ is an affine open embedding. Then $Rj_!$ and $Rj_*$ are both perverse $t$-exact functors for $\bDbc(X,\bC)$. Consider the  distinguished triangle in $\bDbc(X,\C)$:
$$Rj_! j^{-1}\sF_m \rightarrow   \sF_m \to i_* i^{-1}\sF_m \overset {[1]}{\to}.$$ 
Due to the choice of the stratification, $j^{-1}\sF_m$ is a shifted local system, perverse sheaf, by perverse $t$-exactness, $Rj_! j^{-1}\sF_m \in \Perv(X, \C)$ for any $m\subset R$. We then have a long exact sequence of perverse sheaves
\begin{center}
$0  \to  { }^p \sH^{-1} (i_* i^{-1}\sF_m)  \to Rj_! j^{-1}\sF_m \to  \sF_m  \to { }^p \sH^0 (i_* i^{-1}\sF_m) \to 0.$
\end{center} 
Hence $$\ell(\sF_m)=\ell(Rj_! j^{-1}\sF_m)+ \ell({ }^p \sH^0 (i_* i^{-1}\sF_m))-\ell({ }^p \sH^{-1} (i_* i^{-1}\sF_m)).$$

Now the supports of $i_*i^{-1}\sF_R$ and $i_* i^{-1}\sF_m$   for any $m\subset R$ are included in $Z$ by the choice of simultaneous stratification. Here $\dim Z<n$. Note that $Rj_*$, $j^{-1}$, $i_*$, $i^{-1}$,  and ${}^p \sH^k$  are (absolute) $\Q$-constructible functors with respect to the obvious unispaces (\cite[Theorem 6.4.3]{BW17}). In particular, they satisfy generic base change for $\otimes^L_RR/m$. Hence 
$$\ell({ }^p \sH^k (i_* i^{-1}\sF_m))=\ell({ }^p \sH^k (i_* i^{-1}\sF_R)\otimes^L_RR/m)$$
and
$$\ell(Rj_! j^{-1}\sF_m)=\ell((Rj_! j^{-1}\sF_R)\otimes^L_RR/m)$$
for $m$ generic in $\spec(R)(\bC)$. Therefore the function
$$
m\mapsto \ell({ }^p \sH^0 (i_* i^{-1}\sF_m))-\ell({ }^p \sH^{-1} (i_* i^{-1}\sF_m))
$$
is $\Q$-constructible in the classical sense on $\spec(R)(\bC)$ by induction, and we are left to showing the same for the function 
$$m\mapsto \ell(Rj_! j^{-1}\sF_m)=\ell((Rj_! j^{-1}\sF_R)\otimes^L_RR/m).$$

Again as above, we can assume that $j^{-1}\sF_R$ is a shifted local system of free finitely generated $R$-modules such that $$j^{-1}\sF_R\otimes^L_RR/m=j^{-1} \sF_m$$ for all $m$. Set $L$ the local system agreeing up to a shift with $j^{-1}\sF_R$, and similarly $L_m$ the one agreeing with $j^{-1} \sF_m$. By the previous argument we can assume that $L$ and all $L_m$ have the same length.
Assume that the local system $L$ is not simple. Equivalently $L_m$ are not simple. Take any non trivial sub-local system $L^\prime_m$ of $L_m$, one have the following short exact sequence:
\begin{center}
$0\to L^\prime_m \to L_m \to L^{\prime\prime}_m \to 0$
\end{center}
Applying the functor $Rj_!$, we get a short exact sequence of perverse sheaves:
\begin{center}
$0\to Rj_!(L^\prime_m [n]) \to Rj_!( L_m[n])=Rj_!j^{-1}\cF_m \to Rj_!( L^{\prime\prime}_m[n]) \to 0$
\end{center}
In particular, $$\ell(Rj_!j^{-1}\sF_m)= \ell(Rj_!(L^\prime_m [n]))+\ell(Rj_!( L^{\prime\prime}_m[n])).$$
Based on this observation, it is easy to see that $$\ell(Rj_!j^{-1}\sF_m)=\ell(Rj_!(q(L_m) [n])),$$
where
$$
q:\sR(X,r)\lra \cM(X,r)
$$
is the surjection to the categorical quotient of $\sR(X,r)$ by the conjugation action of $GL_r(\bC)$. Indeed, it is well-known that $\cM(X,r)$ is a coarse moduli space of local systems whose $\bC$-points parametrize the isomorphism classes of semi-simple $\bC$-local systems on $X$. Moreover, $q$ associates to a representation its semi-simplicization. $\cM(X,r)$ is an affine scheme defined over $\bQ$ and $q$ is a morphism of schemes defined over $\bQ$.

By the harmless assumptions made above,
$$\ell(Rj_!j^{-1}\sF_m)=\ell(Rj_!(q(L) [n])\otimes^L_RR/m),$$
at least generically for $m\in \spec(R)(\bC)$, which is all that matters anyway.

Now consider the distinguished triangle in $\bDbc(X,R)$ 
\begin{center}
$Rj_!(q(L)[n]) \to j_{!*}(q(L)[n]) \to \sK\xra{[1]}$
\end{center}
where the functor intermediate extension $j_{!*}$ is defined for $R$-coefficients  inductively from pullbacks and direct images along strata together with trunctations and shifts (recall that we have used a similar definition above for defining perverse cohomology functors for $R$-coefficients). Here $\sK$ is the cone and has support in $Z$. Since $j_{!*}$ is an (absolute) $\bQ$-constructible functor, see {\it loc. cit.}, 
$$
j_{!*}(q(L)[n])\otimes^L_RR/m=j_{!*}(q(L_m)[n])
$$
for $m$ generic. In particular,
$$
\ell(j_{!*}(q(L)[n])\otimes^L_RR/m=\ell(j_{!*}(q(L_m)[n]))
$$
for $m$ generic. Since $q(L)$ is a semi-simple local  system and $j_{!*}(q(L_m)[n])$ is the intersection complex of $q(L_m)$,
$$\ell(j_{!*}(q(L_m)[n]))= \ell(q(L_m))=\ell(L_m).$$ So $$\ell(Rj_!(q(L) [n])\otimes^L_RR/m)=\ell(j_!(L_m [1]))= \ell(L_m)+\ell(\sK\otimes^L_RR/m),$$
at least for generic $m$.

The function $m\mapsto \ell(\sK\otimes^L_RR/m)$ is $\bQ$-constructible on $\spec(R)(\bC)$ by induction, since the dimension of the support of $\sK$ is less than $n$. 
 The function $m\mapsto \ell(L_m)=\ell(L\otimes_RR/m)$ was shown to be $\bQ$-constructible earlier in this proof. It follows that the function
$$
m\mapsto \ell(Rj_!(q(L) [n]) 
$$
is $\bQ$-constructible and this finishes the proof of the theorem.

\section{Example}

Currently one does not have  a nice interpretation of absolute $\bQ$-constructible sets of semi-simple local systems of rank higher than 1. We provide though here some examples of rank 2 loci with prescribed length for the associated direct images as in Corollary \ref{corJ}.

Let $$j:U=\C \setminus \{ a, b \}\ra\C$$ be the inclusion, where $a$ and $b$ are two different points in $\C$. The fundamental group of  $U$ is  the free group with two generators.
We focus on the case of representations of $\pi_1(U)$ into $SL_2(\C)$. The moduli space of semi-simple rank 2 local systems, denoted $$\sM(U, SL_2),$$ can be then computed to be  $\C^3$ with the coordinates $(x,y,z)$ (\cite[Theorem A]{Gol}), where
\begin{center}
$x= tr(A), y=tr(B), z=tr(AB)$.
\end{center}
Here  $A$ and $B$ denote the  $2\times 2$ matrices associated to the generators given by loops around the points $a$ and $b$, respectively.

The simple $SL_2(\bC)$-local systems are those in the complement of  the zero locus of the polynomial
\begin{center}
$z^2-xyz-4+y^2+x^2$.
\end{center}

One has a short exact sequence of perverse sheaves for a $\bC$-local system $L \in \sM(U, SL_2)$:
\begin{center}
$0 \to j_{!*}(L[1]) \to Rj_*(L [1]) \to \sQ \to 0$
\end{center}
where $\sQ= \sQ_a \oplus \sQ_b$. Here $\sQ_a $ and $ \sQ_b$ are skyscraper sheaves supported on the points $a$ and $b$ with stalk isomorphic to $H^1(U_a, L)$ and $H^1(U_b, L)$, respectively, where $U_a$ and $U_b$ are small punctured balls centered at $a$ and $b$. Note that  $L$ is a semi-simple local system, hence $\ell(j_{!*}(L [1]))= \ell(L )$.  So \begin{center}
$\ell(Rj_*(L [1]) ) =\ell(L) + \dim H^1(U_a, L) + \dim H^1(U_b, L).$
\end{center}
 Moreover   $H^1(U_a, L)$ and $H^1(U_b, L)$ coincide with the number of Jordan blocks with eigenvalue 1 for $A$ and $B$, respectively.
 
\medskip

Now one can compute the list of the loci of local systems $L$ in $\sM(U, SL_2)\simeq\bC^3$ with prescribed length for  the perverse sheaves $Rj_*(L [1])$:

\begin{itemize}
\item $ \{ L \in \sM(U, SL_2) \mid \ell(Rj_*(L [1])) \geq 7 \}=\emptyset.$
  
\item $ \{ L \in \sM(U, SL_2) \mid \ell(Rj_*(L [1])) \geq 6 \} =\{(2,2,2)\}. $  

\item $\{ L \in \sM(U, SL_2) \mid \ell(Rj_*(L [1])) \geq 5 \} =\{(2,2,2)\}.$
  
  \item $\{ L \in \sM(U, SL_2) \mid \ell(Rj_*(L [1])) \geq 4 \} = \{x=2, y=z\} \cup \{y=2, x=z\}. $
  
  \item $\{ L \in \sM(U, SL_2) \mid \ell(Rj_*(L [1])) \geq 3 \} =\{ x=2, y=z \} \cup \{x=2, y=2\} \cup \{y=2, x=z\}. $
  
  \item $\{ L \in \sM(U, SL_2) \mid \ell(Rj_*(L [1])) \geq 2 \} =\{ z^2-xyz+x^2+y^2-4=0 \} \cup \{x=2\} \cup \{y=2\}. $
\end{itemize} 
It is interesting to observe that, on the simple locus, these length jump loci of rank 2 local systems are linear.


\end{document}